
\baselineskip=14pt
\parskip=10pt

\magnification=\magstephalf

\def\1{{\overline{1}}}
\def\2{{\overline{2}}}
\parindent=0pt
\overfullrule=0in

\def\frac#1#2{{#1 \over #2}}
\centerline
{\bf 
Automated Counting and Statistical Analysis of Labeled Trees with Degree Restrictions
}
\bigskip
\centerline
{\it Shalosh B. EKHAD and Doron ZEILBERGER}
\bigskip

{\bf Maple Package}

This article is accompanied by a Maple package, {\tt ETSIM.txt}, available from

{\tt https://sites.math.rutgers.edu/\~{}zeilberg/tokhniot/ETSIM.txt} \quad .

The web-page of this article,

{\tt https://sites.math.rutgers.edu/\~{}zeilberg/mamarim/mamarimhtml/etsim.html} \quad ,

contains input and output files, referred to in this paper.

{\bf Sneak Preview}

In 1887 Arthur Cayley famously proved that there are $n^{n-2}$ labeled trees on $n$ vertices, but how many trees are there where, for example, all the
vertex-degrees {\bf must} belong to the set $\{1,3,4,6,7\}$? Typing

{\tt LTseq($\{$1,3,4,6,7 $\}$,50);}

in the Maple package {\tt ETSIM.txt} accompanying the present article, you would get the first 50 terms immediately.
Typing

{\tt LTseq($\{$1,3,4,6,7 $\}$,1000)[1000];}

would give you, just as fast, the $2784$-digit integer that gives you the {\bf exact} number of labeled trees with $1000$ vertices such that every vertex either has one neighbor, 
or three neighbors, or four neighbors,  or six neighbors, or seven neighbors.

On the other hand, if you want the first $70$ terms of the sequence enumerating, for example, labeled trees where
{\bf none} of the vertices has degrees in the set $\{2,3,5\}$,  just type:

{\tt LTseqF($\{2$,3,5$\}$,70);} \quad.

{\bf An Important Formula}

Using the general {\it generatingfunctiology} for {\bf labeled objects}, and using {\bf Lagrange Inversion}  (see [Z] for a lucid and engaging account), one
can easily establish the following theorem.

{\bf Important Theorem}: Let $P$ be a (finite or infinite) set of positive integers. Let $a_P(n)$ be the number of labeled trees where each vertex must have a number of
neighbors that belongs to $P$, then

$a_P(n)$= (n-2)! $\cdot$ Coefficient of $z^{n-2}$ in $\left (  \sum_{i \in P} \frac{z^{i-1}}{(i-1)!}     \right )^n $ \quad .

{\bf Part I: Efficient Counting, Recurrences, and Asymptotics for Enumerating Labeled trees where all vertices MUST belong to a given  Finite Set}

In this case the amazing {\bf Almkvist-Zeilberger} algorithm [AZ] (see [D] for a lucid and engaging exposition) can easily find a linear recurrence equation with polynomial coefficients
satisfied by $a_P(n)$. This is implemented in procedure {\tt LTrec} in the Maple package. This also enables one to find asymptotics. Procedure {\tt InfoLT} 
takes care of this. It also gives you the limiting distribution of the respective participating degrees.

To get an article with information about {\bf all}  possible subsets of $\{1,2,\dots, M\}$ that include $1$ (of course), with the first $K1$ terms 
of each sequence displayed, type 

{\tt LTpaper(M,K1,n);} \quad .

Here is the output file when you choose {\tt M=7} and {\tt K1=30}:

{\tt https://sites.math.rutgers.edu/\~{}zeilberg/tokhniot/oETSIM1a.txt} \quad.

The output file

{\tt https://sites.math.rutgers.edu/\~{}zeilberg/tokhniot/oETSIM1.txt} \quad

only does it for subsets of $\{1,2,3,4,5\}$, but also gives the limiting distribution of the respective participating degrees, that is much more time-consuming.

{\bf Part II: Efficient Counting, and Estimated Asymptotics, for Enumerating Labeled trees where all vertices must NOT belong to a  specified Finite Set}

If $F$ is the forbidden set, then in the Important Formula, $P=\{1,2,3, \dots \} \backslash F$, and we have that the number of labeled trees on $n$ vertices
such that none of the vertices have a degree in $F$ is

(n-2)! $\cdot$ Coefficient of $z^{n-2}$ in $\left (  e^z - \sum_{i \in F} \frac{z^{i-1}}{(i-1)!}     \right )^n $ \quad .

If $F$  is a singleton $\{ r \}$ we have that the the number of labeled trees on $n$ vertices such that no vertex has $r$ neighbors is

(n-2)! $\cdot$ Coefficient of $z^{n-2}$ in $\left (  e^z - \frac{z^{r-1}}{(r-1)!}     \right )^n $ \quad .

This equals

(n-2)! $\cdot$ Coefficient of $z^0$ in $\left (  e^z - \frac{z^{r}}{(r-1)!}     \right )^n \cdot \frac{1}{z^{n-2}}$ \quad ,

that equals

(n-2)! $\cdot$ Coefficient of $z^0$ in $ \sum_{k=0}^{n} (-1)^k \, {{n} \choose {k}} \, \frac{1}{(r-1)!^k} \, \frac{e^{z(n-k)}}{z^{n-2-r\,k}} $ \quad ,

that equals

$$
\sum_{k=0}^{\lfloor (n-2)/r \rfloor} \, (-1)^k  {{n} \choose {k}} \, \frac{1}{(r-1)!^k} \, \frac{(n-2)! (n-k)^{n-2-r\,k} }{(n-2-rk)!} \quad ,
$$
that enables a very fast computation. This is implemented in procedure {\tt T1rn(r,n)}. If the forbidden set has two elements, we have a double sum,
implemented in procedure  {\tt T2rn(r1,r2,n)}, and if $F$ has three elements then we have a triple sum (procedure {\tt T3rn(r1,r2,r3,n)}).

Procedure {\tt LTseqF(F,N)} gives the first {\tt N} terms of the enumerating sequence of labeled trees avoiding the members of $F$ as vertex-degrees.

Procedure {\tt LTFpaper} gives information about many cases, together with {\it estimated} (non-rigorous, but very reliable!) asymptotics.

For many such sequences, and estimated asymptotics, see

{\tt https://sites.math.rutgers.edu/\~{}zeilberg/tokhniot/oETSIM2.txt} \quad .

{\bf Part III: Statistical Analysis}

Given a labeled tree $T$, Let $X_d(T)$  denote the number of vertices that have degree $d$. It is a {\it random variable},  defined on the {\it sample space} of 
all the $n^{n-2}$ labeled trees on $n$ vertices.
In particular $X_1(T)$ is the number of leaves of $T$. 

We are interested in the {\it expectation}, {\it variance}, and higher moments of $X_d$. We will show that for each $d$, $X_d$ is {\it asymptotically normal}.

We will also explore  how $X_{d_1}$ and $X_{d_2}$ interact, as $n$ goes to infinity.

These questions are answered by the following very interesting theorem.

{\bf Interesting Theorem:} Let $d$ be a positive integer, and let $X_d(T)$ be the number of vertices of a labeled tree $T$, that have degree $d$, then
we have the following interesting facts.

$\bullet$ 
$$
E[X_d]  \, = \, \frac{e^{-1}}{(d-1)!} \cdot n + O(1)  \quad .
$$

$\bullet$ 
$$
Var(X_d)= \left ( \frac{e^{-1}}{(d-1)!} - \frac{( d^2-4d+5)e^{-2}}{(d-1)!^2} \right ) n + O(1) \quad .
$$

$\bullet$ $X_d$ is {\it asymptotically normal}.

$\bullet$ If $1 < d_1 < d_2$, then 
$$
Cov(X_{d_1},X_{d_2}) \, = \,  \frac{2d_1+2d_2-d_1d_2-5}{(d_1-1)! (d_2-1)!e^2} \cdot \, n + O(1) \quad .
$$

$\bullet$ $X_{d_1}$ and $X_{d_2}$ are joint-asymptotically normal with  (limiting) correlation coefficient
$$
\frac{
 2d_1+2d_2-d_1d_2-5
}
{
\sqrt
{
 (\sqrt{(d_1-1)!} - ( d_1^2-4d_1+5)e^{-1})
\,
 (\sqrt{(d_2-1)!} - ( d_2^2-4d_2+5)e^{-1})
}
} \quad .
$$

{\bf Comment added Feb. 1, 2022}: We found out that the expressions for the expectation and variance of $X_d$, as well as the fact that it is (singly-) asymptotically normal, go back to Alfr\'ed R\'eny [R], for the $d=1$ case,
and to Amram Meir and John W. Moon [MM] for the general case. See also theorem 7.7 (p. 73) of Moon's classic monograph [M]. We hope that the result about the
covariance, and the {\it joint} asymptotic normality of $X_{d_1}$ and $X_{d_2}$ are new. We would appreciate any references, of course, in case our hope is wrong.

{\bf Sketch of Proof:}

Define the {\bf weight} of a tree $T$ to be 
$$
wt(T):=\prod_{v \in T} g_{degree(v)} \quad,
$$ 
where $g_1, g_2, \dots$ are {\bf commuting indeterminates}. The same reasoning that lead to the
Important Formula tells you that the weight-enumerator of labeled trees with $n$ vertices is

$(n-2)! \cdot$ Coefficient of $z^{n-2}$ in $\left (  \sum_{i=1}^{\infty}  \frac{g_i z^{i-1}}{(i-1)!}     \right )^n $ \quad .

Since we want to focus on $X_d$, we set all the $g_i$'s to $1$ except for the {\it active} variable, $g_d$.
This gives us that the {\it probability generating function} (under the uniform distribution on labeled trees) of $X_d$, let's call it
$p_d(g_d)$, is:

$\frac{(n-2)!}{n^{n-2}}$ $\cdot$ Coefficient of $z^{n-2}$ in $\left (e^{z}+ (g_d-1)\frac{z^{d-1}}{(d-1)!}    \right )^n $ \quad .

The expectation, $E[X_d]$, is $\frac{d}{dg_d} p(g_d) \Bigl \vert_{g_d=1}$, and more generally, the $k$-th moment, $E[X_d^k]$ is 
 $(g_d \frac{d}{dg_d})^k p(g_d) \Bigl \vert_{g_d=1}$. For any {\it specific} $d$ and $k$, Maple can find explicit expressions, as finite linear combinations of
terms of the form $(n-i)^{n-j}/n^{n-2}$, from which Maple can find, in turn,  explicit expressions for the {\it moments about the mean} for
each {\it numeric} $d$ and each {\it numeric} $k$. From these, for small $k$, one can easily guess explicit expressions for
general $d$, and with a little more effort, one can even get the computer to do it for {\bf symbolic} $d$, but still numeric (small) $k$. 
As $k$ gets larger, the expressions get more and more complicated, but for deducing {\it limiting behavior}, the leading terms suffice,
and it is conceivable that one may be able to do it for {\bf both} symbolic $d$ and symbolic $k$, but we decided that
we have better things to do.

Similarly for mixed moments. See procedure

{\tt LTmom2am(d1,d2,n,k1,k2)}  \quad,

that gives you an {\bf explicit} expression, in $n$, for the mixed $(k_1,k_2)$ moment, about the mean, of $X_{d_1}$ and $X_{d_2}$, but as $k_1$ and $k_2$
get larger, the expressions become larger and larger. On the other hand, Procedure

{\tt LTmom2amL(d1,d2,k1,k2)}  \quad ,

gives us the {\bf exact} value of the limit, as $n$ goes to infinity, of the former divided by $n^{ \lfloor (d_1+d_2)/2  \rfloor }$, from which it is easy
to guess {\it symbolic} expressions in $d_1,d_2$ for each specific, numeric, $k_1$ and $k_2$.

Since we are  {\bf experimental mathematicians}, doing it for sufficiently many moments, and comparing the limiting scaled moments to those of the
normal distribution, and the limiting scaled mixed moments to the mixed moments of a bivariate normal with the above limiting correlation, is good enough for us!

For lots of nice formulas, read the output file

{\tt https://sites.math.rutgers.edu/\~{}zeilberg/tokhniot/oETSIM3.txt} \quad .

Enjoy!

{\bf References}

[AZ]  Gert Almkvist and Doron Zeilberger, {\it The method of differentiating Under The integral sign}, J. Symbolic Computation {\bf 10} (1990), 571-591. \hfill \break
{\tt https://sites.math.rutgers.edu/\~{}zeilberg/mamarim/mamarimPDF/duis.pdf} \quad .

[D]  Robert Dougherty-Bliss,  {\it Integral recurrences from A to Z}, Amer. Math. Monthly, to appear. {\tt https://arxiv.org/abs/2102.10170} \quad .

[MM] A. Meir and J.W. Moon, {\it On nodes of degree two in random trees}, Mathematika {\bf 15}, 188-192.

[M] J.W. Moon, {\it ``Counting Labelled Trees''}, Candadian Mathematical monographs \#1, Canadian Mathematical Congress, 1970.

[R] A. R\'enyi, {\it Some remarks on the theory of trees}, Publications of the Mathematical Institute of the Hungarian Academy of Science {\bf 4} (1959), 73-85.

[Z] Doron Zeilberger, {\it Lagrange Inversion Without Tears (Analysis) (based on Henrici)}, Personal Journal of Shalosh B. Ekhad and
Doron Zeilberger, 2002.\hfill\break
{\tt https://sites.math.rutgers.edu/\~{}zeilberg/mamarim/mamarimhtml/lag.html} \quad .

\bigskip
\hrule
\bigskip
Shalosh B. Ekhad and Doron Zeilberger, Department of Mathematics, Rutgers University (New Brunswick), Hill Center-Busch Campus, 110 Frelinghuysen
Rd., Piscataway, NJ 08854-8019, USA. \hfill\break
Email: {\tt [ShaloshBEkhad, DoronZeil] at gmail dot com}   \quad .

{\bf Exclusively published in the Personal Journal of Shalosh B. Ekhad and Doron Zeilberger and arxiv.org}

First Published: {\bf  15 of Shevat, 5782} (Birthday for the Trees; Jan. 17, 2022, according to their reckoning).  This version: {\bf Feb. 1, 2022.}

\end